\documentclass{amsart}
\usepackage[dvips]{graphicx}
\input xy
\xyoption{all}

\newtheorem{theorem}{Theorem}
\newtheorem{lemma}[theorem]{Lemma}
\newtheorem{proposition}[theorem]{Proposition}

\theoremstyle{definition}

\numberwithin{equation}{section}

\newcommand{\zz}{\mathbb{Z}}

\newcommand{\hnu}{{\sf n}}
\newcommand{\mk}{{\sf m}(K)}
\newcommand{\lk}{{\sf l}(K)}
\newcommand{\sig}{{\sf S}_{K}}

\begin{document}

\title[Homologous Non-isotopic Symplectic
Surfaces of Higher Genus]{Homologous Non-isotopic Symplectic\\
Surfaces of Higher Genus}

\author{B. Doug Park}
\address{Department of Pure Mathematics, University of Waterloo, Waterloo,
Ontario, N2L 3G1, Canada} \email{bdpark@math.uwaterloo.ca}
\thanks{B. D. Park was partially supported by NSERC and CFI/OIT grants.}

\author{Mainak Poddar}
\address{Department of Pure Mathematics, University of Waterloo, Waterloo,
Ontario, N2L 3G1, Canada} \email{mpoddar@math.uwaterloo.ca}

\author{Stefano Vidussi}
\address{Department of Mathematics, University of California,
Riverside, CA 92521, USA}
\email{svidussi@math.ucr.edu} \thanks{S. Vidussi was partially
supported by NSF grant \#0306074.}
\dedicatory{Dedicated to Ron Fintushel on the occasion of his sixtieth birthday}
\subjclass[2000]{Primary 57R17, 57M05; Secondary 53D35, 57R95}
\date{Revised on July 22, 2005}

\begin{abstract}
We construct an infinite family of homologous, non-isotopic,
symplectic surfaces of any genus greater than one in a certain
class of closed, simply connected, symplectic four-manifolds.  Our
construction is the first example of this phenomenon for surfaces
of genus greater than one.
\end{abstract}

\maketitle

\section{Introduction}\label{sec:introduction}

Several papers have addressed, in the last few years, the problem
of isotopy of symplectic surfaces representing the same homology
class $\alpha \in H_{2}(M;\zz)$ of a simply connected, symplectic
$4$-manifold $M$. Although a list of the results obtained in those
papers is beyond the scope of this introduction, we just want to
mention the existence of both uniqueness and non-uniqueness
results, depending on the choice of the pair $(M,\alpha)$. The
prototypical example of a non-uniqueness result, further discussed
in this paper, is the presence of infinitely many non-isotopic
symplectic tori representing any multiple of the fiber class of
(some) elliptic surfaces. The idea behind this construction,
developed by Fintushel and Stern in \cite{fs:non-isotopic}, is to
obtain nonequivalent symplectic tori by braiding parallel copies
of a symplectic torus.  All known examples of non-uniqueness, so
far, have been genus $1$ representatives of a homology class of
self-intersection $0$. On the other hand, the few higher genus
cases under control have always led to uniqueness results as in
\cite{ST}. The interest in obtaining examples of non-isotopic
symplectic surfaces of higher genus has been pointed out by many
researchers, in particular, in \cite{ADK}, \cite{fs:non-isotopic}
and \cite{S2}.

There is no special reason to expect that non-isotopy phenomena
should be restricted to genus $1$, self-intersection $0$ surfaces.
In fact, if we also take into account $4$-manifolds with
non-trivial fundamental groups, Smith has shown in \cite{smith}
the existence of higher (odd) genus examples, and examples of
symplectic curves with cusp singularities have been provided (for
the projective plane) in \cite{ADK}. However, all attempts to
build higher genus examples in simply connected $4$-manifolds by
suitably ``doubling" the braiding construction have failed. The
reason for the failure is the impossibility of detecting
non-isotopy (for potential examples) either by using the
Seiberg-Witten invariants or by using a more classical topological
invariant like the fundamental group of the complement. In the
latter approach, the reason is that most of the basic building
blocks, genus $1$ surfaces built through the braiding
construction, cannot be distinguished by means of the $\pi_{1}$
either, as first observed in \cite{fs:non-isotopic}.

However, there is a second mechanism, different from the original
braiding construction, that allows us to produce non-isotopic
symplectic tori in a symplectic $4$-manifold (see \cite{ep:E(1)_K}
and \cite{V2}).  This alternative construction originates from the
presence of non-isotopic, nullhomologous Lagrangian tori in a
large class of (simply connected) symplectic $4$-manifolds (see
\cite{FS3} and \cite{V1}). These symplectic tori are obtained by
``summing" a preferred symplectic torus with non-isotopic
nullhomologous tori. The non-isotopy can be detected, as for the
braiding construction, by using Seiberg-Witten theory, but the
advantage of this second construction is that, in many cases,
$\pi_{1}$ is sufficient to distinguish the tori, as shown in
\cite{ep:fundamental group}.  And if we ``double" the
construction, the fundamental groups of the complements may retain
enough information to distinguish the resulting higher genus
surfaces.

In fact, following the ideas outlined above, we will be able to
prove the following.

\begin{theorem}\label{main}
For any fixed integer\/ $q \geq 1$, there exist simply connected,
symplectic $4$-manifolds $X$ containing infinitely many
homologous, pairwise non-isotopic, symplectic surfaces\/
$\{\Xi_{p,q} \mid \gcd(p,q)=1\}$ of genus\/ $q+1$. Furthermore,
there is no pair homeomorphism between $(X,\Xi_{p,q})$ and
$(X,\Xi_{p',q})$ unless $p'=p$.
\end{theorem}

We briefly preview the idea of the construction. First, consider
the symplectic simply connected 4-manifold $E(2)_{K}$ obtained by
the knot surgery construction of \cite{fs:knots} from the elliptic
$K3$ surface $E(2)$ and a non-trivial fibered knot $K$ in $S^3$.
$E(2)_{K}$ contains, for each $q \geq 1$, an infinite family of
homologous, symplectic, non-isotopic tori $\{T_{p,q}  \mid
\gcd(p,q)=1\}$, together with a symplectic surface $\Sigma_{g}$ of
genus $g = g(K) + 1$ of self-intersection $0$. This surface
$\Sigma_{g}$ intersects each of the tori $T_{p,q}$ at $q$\/
positive transverse points. By doubling $E(2)_{K}$ along
$\Sigma_{g}$ we obtain a symplectic, simply connected $4$-manifold
$X=D_{K}$. By summing together two copies of the torus $T_{p,q}$,
one copy from each side, we obtain a family of homologous genus
$q+1$, self-intersection $0$ symplectic surfaces $\Xi_{p,q}$.

By direct computation, we will show that if the fibered knot $K$\/
used in the construction of $D_{K}$ is non-trivial, then
infinitely many of these surfaces $\Xi_{p,q}$ have
non-homeomorphic complements  $D_K\!\setminus\Xi_{p,q}$,
distinguished by their fundamental groups. This proves, in
particular, that the surfaces are not isotopic.

\section{Construction}\label{sec:construction}

In this section we will construct an infinite family of
non-isotopic, symplectic surfaces of genus greater than or equal
to $2$\/ for a class of simply connected, symplectic
$4$-manifolds. These 4-manifolds are the manifolds denoted by
$D_{K}$, where $K$ is any fibered knot, introduced in
\cite{park:double}.

First we recall the construction of $D_K$.  Let $K\subset S^3$ be
a fibered knot and let $g(K)$ denote its genus.  Let $\Sigma_{K}$
denote the fiber of the fibration of the knot exterior, the
minimal genus Seifert surface for the knot $K$.

For each fibered knot $K$ we can construct a symplectic
$4$-manifold $E(2)_{K}$, homeomorphic to the elliptic surface
$E(2)$ (the $K3$ surface), obtained by knot surgery on $E(2)$ with
the knot $K$ (see \cite{fs:knots}). Denote by $N_{K}$ the fibered
$3$-manifold obtained by the $0$-surgery on $S^{3}$ along $K$.
$S^{1} \times N_{K}$ is a symplectic $4$-manifold with a framed,
self-intersection $0$ symplectic torus $S^{1} \times {\bf m}$,
where ${\bf m}$ is the core of the surgery solid torus. The
manifold $E(2)_{K}$ can be presented as
\begin{equation} \label{knotsur}
E(2)_{K} \:=\: \big(E(2) \setminus \nu F\big) \cup \big(S^{1}
\times (S^{3} \setminus \nu K)\big) \:=\: E(2) \#_{F = S^{1}
\times {\bf m}} S^{1} \times N_{K}\, ,
\end{equation}
where the gluing diffeomorphism identifies factorwise the boundary
3-tori $F\times
\partial D^{2}$ and $S^{1} \times \mu(K) \times \lambda(K)$
(reversing the orientation on the last factor). The presentation
of $E(2)_{K}$ as a fiber sum shows, by Gompf's theory (see
\cite{gompf:sum}), that it admits a symplectic structure
restricting, outside the gluing locus, to the one of the summands.

Inside $E(2)_K$ a disk section $S$\/ of $(E(2) \setminus \nu F)$
can be glued to a Seifert surface $\Sigma_K$ to form a closed,
symplectic, genus $g(K)$ surface $\Sigma$ of self-intersection
$-2$.  By taking $\Sigma$ and a regular torus fiber $F$\/ and
resolving their normal intersection, we obtain a symplectic
surface $\Sigma_g$ of self-intersection $0$ and genus $g=g(K)+1$
inside $E(2)_K$. Such a surface can be endowed with a natural framing, inherited from the canonical framings of the fibers of $E(2)$ and $S^{1} \times N_{K}$. We define
\begin{equation}\label{D_K construction}
D_K \:=\: \big(E(2)_K \setminus \nu \Sigma_g \big) \cup_{\varphi}
\big( E(2)_K \setminus \nu \Sigma_g \big),
\end{equation}
where the gluing map $\varphi$ is an orientation-reversing
self-diffeomorphism on $\partial ({\nu \Sigma_g}) \cong \Sigma_g
\!\times \partial D^{2}$ that is the identity on the $\Sigma_g$ factor and
complex conjugation on the $\partial D^{2}$ factor.

The following result is proved in \cite{park:double}.

\begin{proposition}
$D_K$ is a closed, simply connected, spin, irreducible,
symplectic\/ $4$-manifold.  The signature of\/ $D_K$ is\/ $-32$,
and its intersection form is given by $\,4E_8 \oplus
(7+2g(K))\left[\begin{smallmatrix}0&1\\1&0
\end{smallmatrix}\right]$.
\end{proposition}

In order to obtain higher genus symplectic surfaces in $D_{K}$, we
will start from some symplectic tori in $E(2)_{K}$. The symplectic
$4$-manifold $E(2)_{K}$ contains, for every $q \geq 1$, an
infinite family $\{T_{p,q} \}$ of pairwise non-isotopic symplectic
tori representing the homology class $q[F] = q[S^{1} \times {\bf
m}] \in H_2(E(2)_K;\zz)$ indexed by the integers $p$ coprime to
$q$ (see \cite{ep:E(1)_K} and \cite{V2}).

There are two equivalent ways to present these tori. We start with
the approach of \cite{V2}. First, by looking at the construction
of $E(2)_{K}$ in (\ref{knotsur}), we can recognize a rim torus
$R$\/ given by the image in $E(2)_{K}$ of $S^{1} \times
\lambda(K)$. As we can assume, up to isotopy, that a copy of the
longitude $\lambda(K)$ lies on the fiber of $N_{K}$ (given by the
Seifert surface $\Sigma_{K}$ capped off with a disk of the Dehn
filling), the torus $R$\/ can be assumed to be Lagrangian. Inside
$S^{1} \times N_{K}$ the torus $S^{1} \times \lambda(K)$ is the
boundary of a solid torus (as $\lambda(K)$ bounds a disk of the
Dehn filling) but after the knot surgery this is not generally
true anymore for $R$ (for non-trivial $K$).  It remains, however,
nullhomologous, as it bounds the $3$-manifold $S^{1} \times
\Sigma_{K}$.

Take now a copy of $\partial \nu K$ pushed slightly inside $S^{3}
\setminus \nu K$, in such a way that it is transverse to the
fibration and intersects each fiber in a curve isotopic to the
longitude. The simple closed curves on this torus are the
$(p,q)$-cables of $K$, where $\gcd(p,q) = 1$ (remember that their
linking number with $K$ is $q$), with $\lambda(K)$ and $\mu(K)$
themselves being the $(1,0)$- and the $(0,1)$-cables respectively.
Up to isotopy we can assume that, with the exception of
$(1,0)$-cable, every $(p,q)$-cable is a curve transverse to the
fibration, and its intersection with a fiber consists of positive
points for $q \geq 1$. Denote by $K_{p,q}$ the $(p,q)$-cable of
$K$, for any $q \geq 1$. We can say that this closed curve is
obtained by {\it circle sum} (in the sense of \cite{FS3}) of $q$
copies of the meridian and $p$ copies of the longitude of $K$.

\begin{figure}[ht]
\begin{center}
\includegraphics[scale=.5]{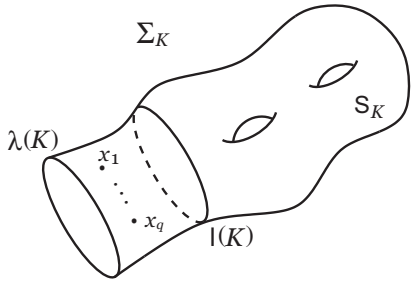}
\end{center}
\caption{$\Sigma_{K}\cap K_{p,q}$ and fibers of the two knot
exteriors} \label{rela}
\end{figure}

It will be useful in what follows to observe that the cable knot
$K_{p,q}$ is contained in an enlarged tubular neighborhood $\hnu
K$ of $K$ (which we can assume to have, say, twice the radius of
$\nu K$ for some metric on $S^{3}$). Note that $S^{3} \setminus
\hnu K$ also fibers over $S^{1}$.  Its fiber $\sig$ is strictly
contained in $\Sigma_{K}$, with $\lk =\partial \sig$ isotopic to
$\lambda(K)=\partial \Sigma_{K}$ inside $\Sigma_{K}$. The cable
$K_{p,q}$ intersects the fiber $\Sigma_{K}$ at $q$ points
$x_{1},\dots,x_{q}$ contained in the knotted annulus $\Sigma_{K}
\setminus {\rm int}(\sig)$, whose boundary is the union
$\lk\cup\lambda(K)$. Figure \ref{rela} illustrates the relations
between these fibers and $K_{p,q}$.

As $K_{p,q}$ is transverse to $\Sigma_{K}$ the torus $S^{1} \times
K_{p,q}$ is symplectic in $S^{1} \times N_{K}$.  For a fixed value
of $q$, the tori $S^{1} \times K_{p,q}$ are homologous, as the
homology class of $K_{p,q}$ in $H_{1}(S^{3} \setminus \nu K;\zz)$
is given by $lk(K_{p,q},K)[\mu(K)] = q [\mu(K)]$ and the Dehn
filling does not affect this relation. (If $q=1$ more is true,
namely the tori $K_{p,1}$ are isotopic, the isotopy being
supported in the solid torus of the Dehn filling.) After knot
surgery, the images $T_{p,q}$ in $E(2)_{K}$ of these tori remain
homologous.  In fact, the image of the homology class of the
$(p,q)$-cable of $K$, through the injective map
\begin{equation*}
H_{1}(S^{3} \setminus \nu K;\zz) \longrightarrow H_{2}(S^{1}
\times (S^{3} \setminus \nu K);\zz) \longrightarrow
H_{2}(E(2)_{K};\zz)
\end{equation*}
is $q[F]$.  However, infinitely many of the tori $T_{p,q}$ become
non-isotopic. Roughly speaking, the torus $T_{p,q}$ is the sum of
$q$ copies of $F$\/ with $p$ copies of the nullhomologous
Lagrangian torus $R$, and can be interpreted in terms of the
cabling of $F$\/ itself.

The family of tori $T_{p,q}$ will be the starting point of our
construction, but note that we can extend the previous argument to
include, instead of just $(p,q)$-cables of the knot $K$, any curve
in $S^{3} \setminus \nu K$ obtained as the circle sum of a number
of copies of the meridian and a closed curve (or a union of closed
curves) lying on the surface $\Sigma_{K}$ (and thus homologous,
but not necessarily isotopic to the longitude). With arguments
similar to the ones above, we can achieve transversality for the
resulting curves; the symplectic tori obtained this way differ by
the addition of nullhomologous Lagrangian tori. Details of this
construction are presented in \cite{V2}.

A second way to construct the family $T_{p,q}$, presented in
\cite{ep:E(1)_K}, is the following. Denote by $H=A\cup B$\/ the
Hopf link. In the knot exterior $S^{3} \setminus \hnu K$ we
consider a standard pair of meridian and longitude, that we
denote, to avoid confusion with $\mu(K)$ and $\lambda(K)$, by
$\mk$ and $\lk$ respectively. We have then the decomposition
\begin{equation}\label{Etgu-Park construction}
E(2)_K  \,=\: \big[E(2)\setminus \nu F \big] \,\cup\,
\big[S^1\!\times(S^3\setminus \nu H)\big]\,\cup\, \big[
S^1\!\times (S^3\setminus \hnu K)\big]\, ,
\end{equation}
with the first gluing identifying $F \times \partial D^2$ with
$S^1 \!\times \mu (A) \times \lambda (A)$ and the second one
identifying $S^1 \!\times \lambda (B) \times \bar{\mu} (B)$ with
$S^1 \!\times \mk \times \lk$ factorwise, reversing the
orientation on the last factor. (The bar over a curve refers to
the negative orientation.) The equivalence between this
construction and the construction of (\ref{knotsur}) (even in the
symplectic category, with a suitable presentation as symplectic
fiber sum) simply follows by observing that the exterior of the
Hopf link is diffeomorphic to $S^{1}\! \times (D^{2} \setminus \nu
\{0\})$.

\begin{figure}[ht]
\begin{center}
\includegraphics[scale=.9]{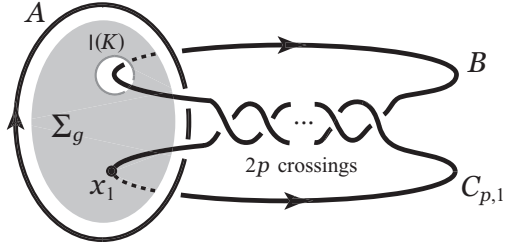}
\end{center}
\caption{3-component link $A\cup B \cup C_{p,1}$ in $S^3$}
\label{fig:doublehopf}
\end{figure}

The torus $T_{p,q}$ is given by the image of $S^1\!\times C_{p,q}
\subset \big[S^1\!\times(S^3\setminus \nu H)\big]$, where
$C_{p,q}$ is the $(p,q)$-cable of $A$, as illustrated in
Figure~\ref{fig:doublehopf}\/ for $q=1$. The transversality of
$C_{p,q}$ to the fibration of $S^{3} \setminus \nu H$ with fiber
given by a disk spanning $A$, pierced once by $\nu B$, is evident,
and this entails that $T_{p,q}$ is symplectic. For fixed $q$, the
curves $C_{p,q}$ are not homologous in $S^{3} \setminus \nu H$,
but they become so after gluing in $S^{3} \setminus \hnu K$ with
the prescription above, as the meridian $\mu(B)$ becomes
nullhomologous. Therefore, the tori $T_{p,q}$ are homologous in
$E(2)_{K}$ for fixed $q$.

The equivalence between the two constructions of $T_{p,q}$ follows
when we observe that the gluing of $S^{3} \setminus \hnu K$ to
$S^{3} \setminus \nu H$ in (\ref{Etgu-Park construction}) is the
{\it splicing}\/ of $K$\/ and $H$\/ along $B$, and the result of
this operation is to produce a copy of $S^{3} \setminus \nu K$,
with the image of $A$\/ giving the knot $K$\/ and the image of the
$(p,q)$-cable of $A$ (i.e. $\!C_{p,q}$) giving the $(p,q)$-cable
of $K$ (i.e. $\!K_{p,q}$).  (See \cite{en}, Proposition 1.1 for a
detailed discussion of this construction.)  Stated otherwise, we
can think of the curve $C_{p,q}$ in the second construction as
spiraling along an inner torus of $S^{3} \setminus \nu H = S^{1}
\times (D^{2} \setminus \nu \{0\})$; after the gluing this inner
torus gets identified with a copy of\/ $\partial \nu K$ pushed
inside $S^{3} \setminus \nu K$, and hence $C_{p,q}$ is identified
with $K_{p,q}$ in the first construction. Note that, after the
splicing, the disk spanning $A$\/ and pierced by $\nu B$ (that
appears shaded in Figure \ref{fig:doublehopf}) is identified with
$\Sigma_{K} \setminus {\rm int}(\sig)$, which appears as the neck
in Figure \ref{rela}.

\begin{figure}[ht]
\begin{center}
\includegraphics[scale=.5]{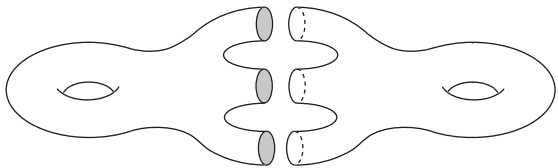}
\end{center}
\caption{Genus 4 surface $\Xi_{p,3}$ as the double of
$(T_{p,3}\setminus\nu\Sigma_g)$} \label{fig:Xi}
\end{figure}

Now that we have defined the family $\{T_{p,q}\mid \gcd(p,q)=1\}$
we can construct, out of two copies of $T_{p,q}$, the desired
genus $q+1$ surface in $D_{K}$. Since in $E(2)_K$,
$[T_{p,q}]\cdot[\Sigma_g]=q[F]\cdot[\Sigma_g]=q[F]\cdot[S] =q\,$
for each $p$, we can internally sum two copies of $q$-times
punctured $T_{p,q}$ from each of $( E(2)_K \setminus \nu \Sigma_g
)$ halves in (\ref{D_K construction}), and thus obtain a family of
homologous, genus $q+1$, self-intersection 0 surfaces
$$\{\,\Xi_{p,q} \:=\: T_{p,q}\# T_{p,q} \:=\:
(T_{p,q} \setminus \nu \Sigma_g ) \cup (T_{p,q} \setminus \nu
\Sigma_g ) \,\}$$ inside $D_K$. (We are going to follow the
convention that $U \setminus V := U \setminus (U \cap V)$.)  See
Figure~\ref{fig:Xi}.

Since each $T_{p,q}$ is symplectic in $E(2)_K$ and our ambient
$4$-manifold $D_K$ is a symplectic sum (cf.$\;$\cite{gompf:sum})
of two copies of $E(2)_K$ along a symplectic surface $\Sigma_g$,
each $\Xi_{p,q}$ is a symplectic submanifold of $D_K$.  In the
next section we discuss how to distinguish the isotopy classes of
submanifolds $\Xi_{p,q}$ by comparing the fundamental groups of
their complements in $D_K$.

\section{Fundamental groups}\label{sec:pi_1}

The isotopy class of the torus $T_{p,q}$ inside $E(2)_{K}$ can be
detected via the Seiberg-Witten invariant of the fiber sum
manifold $E(1)\#_{F=T_{p,q}} E(2)_{K}$ (that depends ultimately on
the diffeomorphism type of the pair $(E(2)_{K},T_{p,q})$).  It
would be interesting (and probably quicker) to be able to use a similar
approach to distinguish the diffeomorphism type of the pair
$(D_{K},\Xi_{p,q})$.  Unfortunately, we currently lack a satisfactory
machinery within the Seiberg-Witten theory to carry this out.
However, we can use a classical approach to this problem, namely the study 
of the fundamental group of the exterior of $\Xi_{p,q}$. This allows us to
prove an even stronger result, namely the existence of pairs $(D_{K},\Xi_{p,q})$ that are not even 
homeomorphic. We were first led to this approach by
the observation that the tori $T_{p,q}$ themselves can be
distinguished by means of the fundamental group alone, as
illustrated in \cite{ep:fundamental group}. The reason underlying
that phenomenon seems to be related to the role of the
nullhomologous Lagrangian tori in the construction of $T_{p,q}$.

It remains an interesting question to determine whether there are homeomorphic
pairs (or at least pairs with the same fundamental group of the exterior) which are
not diffeomorphic. The main result of this section is the following theorem.

\begin{theorem}\label{thm:fundamental}
Let\/ $\Xi_{p,q}$ be the genus\/ $q+1$, symplectic surface inside
$D_K$ constructed in Section~$\ref{sec:construction}$.  Then we
have
\begin{equation}\label{pi_1(Xi_q complement)}
\pi_1(D_K\setminus \Xi_{p,q})\:\cong\: \pi_{1}(S^{3}_{q/p}(K)),
\end{equation}
where $S^{3}_{q/p}(K)$ is the rational homology sphere obtained by
Dehn $(q/p)$-surgery on $S^{3}$ along the knot $K$, so that
\begin{equation} \label{stand} \pi_{1}(S^{3}_{q/p}(K)) \:\cong\:
\frac{\pi_1 (S^3 \setminus \nu K)} {\mu(K)^{q}\lambda(K)^p = 1}
\end{equation} with $\mu(K)$ and\/ $\lambda(K)$ denoting the meridian
and the longitude of\hspace{1pt} $K$ respectively.
\end{theorem}

Theorem \ref{main} now follows from
Theorem~\ref{thm:fundamental}\/ and the fact that when $K$\/ is a
non-trivial knot, for a fixed value of $q$, infinitely many
choices of $p$\/ give mutually non-isomorphic groups
$\pi_{1}(S^{3}_{q/p}(K))$.  A proof of this latter statement is
contained in \cite{ep:fundamental group}. For some knot $K$, we can 
explicitly determine infinite values of $p$ that lead to pairwise non-homeomorphic pairs.
As an example, consider the figure eight knot. This knot is known to be fibered and hyperbolic. 
Moreover, it is known (see \cite{T}) that, with $6$ explicit exceptions all having $|p| < 2$, its $(q/p)$-Dehn surgery is a hyperbolic manifold. At this point we can invoke the result of \cite{BHW} which guarantees that, restricting ourselves to
positive $p$, no two surgeries give homeomorphic $3$-manifolds. Mostow rigidity then implies that all the resulting manifolds have distinct fundamental groups.

The proof of Theorem~\ref{thm:fundamental} will occupy the
remainder of this section. We start by observing that we have the
following decomposition into halves
\begin{equation}\label{decomposition:Xi_q complement}
D_K\setminus\nu \Xi_{p,q} \:\cong\: \big[E(2)_K \setminus
\nu(\Sigma_g\cup T_{p,q}) \big]\:\cup\:\big[E(2)_K \setminus
\nu(\Sigma_g\cup T_{p,q}) \big].
\end{equation}
We want to apply Seifert-Van Kampen theorem to this decomposition.
Our first step is the computation of the fundamental group of the
halves. We claim that
\begin{equation}\label{pi_1(half)}
\pi_1\big(E(2)_K\setminus \nu(\Sigma_g \cup T_{p,q})\big)\:\cong\:
\pi_{1}(S^{3}_{q/p}(K)).
\end{equation}
Recall from (\ref{knotsur}) and the construction in Section
\ref{sec:construction} that we have
\begin{equation}\label{decomposition:T_q complement}
E(2)_K \setminus\nu T_{p,q} \,=\: \big[E(2)\setminus \nu F \big]
\,\cup\, \big[S^1\!\times(S^3\setminus \nu (K \cup K_{p,q}))\big].
\end{equation}

Our strategy is to use decomposition (\ref{decomposition:T_q
complement}) to compute the fundamental group in
(\ref{pi_1(half)}).  Note that the surface $\Sigma_g$ intersects
both components of decomposition (\ref{decomposition:T_q
complement}). We proceed to compute the fundamental group of the
complement of $\nu \Sigma_g$ in each piece separately.  Denote by
$X_{1}$ and $X_{2}$ the two components in decomposition
(\ref{decomposition:T_q complement}). The first result we need is
a straightforward modification of Proposition~2.1 of \cite{fsnon},
and we present its proof here for completeness. Observe that 
$X_1\setminus  \nu \Sigma_g$ (i.e. $X_1 \setminus  (X_1 \cap \nu \Sigma_g)$),
from the definition in (\ref{decomposition:T_q
complement}), does not depend on $K$.

\begin{lemma}\label{lemma:X_1}
Let $\,X_1 = E(2)\setminus \nu F$.  Then $\,\pi_1 (X_1\setminus
\nu \Sigma_g )=1$.
\end{lemma}

\begin{proof}
Note that $X_1\setminus  \nu \Sigma_g$ 
is $[E(2)\setminus\nu F]\setminus \nu \Sigma'$, where $\Sigma'$ is a genus $1$,
self-intersection $0$ surface with boundary, gotten from the union
$F'\cup S$ of another copy of regular fiber and a disk section of
$E(2) \setminus \nu F$ by resolving their normal intersection
point. It is well known that $[E(2)\setminus\nu F]$ is simply
connected (see e.g. $\!$\cite{gompf-stipsicz}).  It now follows
from Seifert-Van Kampen theorem that the fundamental group in
question is normally generated by the meridian of $\Sigma'$.  By
perturbing the elliptic fibration of $E(2)$ if necessary, we may
assume that there is at least one singular cusp fiber which is
topologically a sphere.  This implies that the meridian of
$\Sigma'$ (which is isotopic to a meridian of $S$) bounds a
topological disk, which is a cusp fiber with a disk removed.
\end{proof}

We now consider the remaining piece in decomposition
(\ref{decomposition:T_q complement}). First, observe that $S^{3}
\setminus \nu(K \cup K_{p,q})$ admits a fibration over $S^{1}$,
with fiber $\Sigma_{K,q}$ given by the spanning surface
$\Sigma_{K}$ of $K$ with $q$\/ disjoint disks\/
$\mathbb{D}_{x_1},\dots ,\mathbb{D}_{x_q}\!$ removed, in
correspondence with the intersection points\/ $x_{1},\dots,x_{q}$
of $\Sigma_{K}$ with $K_{p,q}$. These $q$\/ disks all lie in the
knotted annulus $\Sigma_{K} \setminus {\rm int}(\sig)$. The
existence of such fibration should be quite apparent from the
constructions in Section \ref{sec:construction}, and comes from
gluing, through the splicing, the fibration of $S^{3} \setminus
\hnu K$ with fiber $\sig$ and the fibration of $S^{3} \setminus
\nu (H \cup C_{p,q})$ induced by the spanning disk of $A$\/ (that
gives the punctured knotted annulus above). A section of this
fibration is given by a copy of $\mu(K)$.

Denote by $a_{i},b_{i}$, $i=1,\dots,g(K)$, the generators of
$\pi_{1}(\Sigma_{K})$ and denote by $c_{j}$, $j=1,\dots,q$, the
homotopy classes of the loops around the disks\/
$\mathbb{D}_{x_j}$. Locating the base point on $\partial
\hspace{1pt}\sig$ we can assume that representatives of the
generators $a_{i}$ and $b_{i}$ are all contained in $\sig \subset
\Sigma_{K}$.  The fundamental group of $\Sigma_{K,q}$ (a free
group on $2g(K)+q$ generators) is generated by these elements.

\begin{lemma}\label{lemma:X_2}
Let $\,X_2 = \big[ S^1\!\times (S^3\setminus \nu (K \cup K_{p,q}))
\big]$. Then $\,\pi_1 (X_2\setminus \nu \Sigma_g)$ is isomorphic
to the group
\begin{equation} \label{does}
\frac{\langle x \rangle \:\ast\: \pi_1(S^3\setminus \nu (K \cup
K_{p,q}))}{\{ [x,a_i]=[x,b_i]=[x,c_{j}]= 1 \mid i=1,\dots,g(K) \
{\rm and}\ j=1,\dots,q\}}\: ,
\end{equation}
where $\langle x \rangle$ is the free group on one generator, and
$\,a_1, b_1, \dots , a_{g(K)},b_{g(K)},c_{1}, \dots, c_{q}$ are
the generators of $\,\pi_1(\Sigma_{K,q})$.
\end{lemma}

\begin{proof}
Note that $X_{2}$ is the total space of a fiber bundle:
$$
\xymatrix@C=26pt{ \Sigma_{K,q} \;\;\ar[r]  &
  \;\;S^1\!\times(S^3\setminus \nu
(K \cup K_{p,q}))
\ar[d]^{\;\Pi\:=\:{\rm id}\times\pi} \\
 & S^1 \!\times \mu(K) }
$$
Suppose $X_2\cap \Sigma_g = \Pi^{-1}(t_0)\,$ for some point
$\,t_0\in S^1\!\times\mu(K)$.  Then we can write $X_2\cap \nu
\Sigma_g=\Pi^{-1}(\mathbb{D}_{t_0})\,$ for some small disk
$\mathbb{D}_{t_0}\subset S^1\!\times\mu(K)$.  It follows that $X_2
\setminus \nu \Sigma_g$\/ is the total space of the restricted
bundle:
$$
\xymatrix@C=26pt{ \Sigma_{K,q} \;\;\ar[r]  &
  \;\;\big[S^1\!\times(S^3\setminus \nu
(K \cup K_{p,q}))\big]\setminus \Pi^{-1}(\mathbb{D}_{t_0})
\ar[d]^{\;\Pi|} \\
 & [S^1 \!\times \mu(K)]\setminus \mathbb{D}_{t_0} }
$$
The new base\/ $[S^1 \!\times \mu(K)]\setminus \mathbb{D}_{t_0}$
is a torus with a disk removed and hence homotopy equivalent to a
wedge of two circles.  The generator $x$\/ corresponding to the
$S^1$ factor in $S^1\!\times\mu(K)$\/ now no longer commutes with
$\mu(K)$, but $x$\/ still commutes with the generators of\/
$\pi_1(\Sigma_{K,q})$\/ since the map $\Pi$ had trivial monodromy
in the $S^1$ direction.  Presentation (\ref{does}) follows
immediately.
\end{proof}

Now we are ready to prove (\ref{pi_1(half)}). The intersection of
$X_{1} \setminus \nu \Sigma_{g}$ and $X_{2} \setminus \nu
\Sigma_{g}$ is given by the 3-torus $\partial \nu F =  X_{1} \cap
X_{2}$ minus its intersection with $\nu \Sigma_{g}$. This
intersection is given by the neighborhood of a copy of the
meridional circle $\mu(F)$. We have in fact $\partial \nu F
\setminus \nu \Sigma_{g} = \mu(F) \times (F \setminus \nu \{{\rm
pt}\})$, so that
$$\pi_{1}\big((X_{1} \setminus \nu \Sigma_{g}) \cap
(X_{2} \setminus \nu \Sigma_{g})\big) \:\:=\:\: \zz[\mu(F)] \oplus
(\zz[\gamma_{1}] \ast \zz[\gamma_{2}]),$$ where $\gamma_{1}$ and
$\gamma_{2}$ form a homotopy basis for $\pi_{1}(F)$.  The
generators $\mu(F)$, $\gamma_{1}$ and $\gamma_{2}$ are identified,
respectively, with the classes of ${\lambda(K)}$, $\mu(K)$ and
$S^{1}$ through the gluing map in (\ref{knotsur}). By using the
fact that $X_1\setminus
 \nu \Sigma_g$\/ is simply connected (Lemma~\ref{lemma:X_1}), we deduce that
\begin{eqnarray*}
\pi_1 \big( E(2)_K \setminus \nu(\Sigma_g \cup T_{p,q}) \big) &=&
\pi_1 \big( (X_1 \setminus \nu \Sigma_g) \,\cup\,  (X_2\setminus
\nu
\Sigma_g)  \big)  \\
&=& \frac{\pi_1 \big(S^3 \setminus \nu (K \cup K_{p,q})
\big)}{\mu(K)=1, \lambda(K)=1}\; ,
\end{eqnarray*}
where the last equality comes from observing that the elements of
$\pi_{1}(X_{2} \setminus \nu \Sigma_{g})$ identified with the
elements of $\pi_{1}((X_{1} \setminus \nu \Sigma_{g}) \cap  (X_{2}
\setminus \nu \Sigma_{g}))$  become trivial.  This accounts for
the nullhomotopy of generators\/ $x=[S^1]$, $\mu(K)$ and
$\lambda(K)$ appearing in group (\ref{does}).

We are left, therefore, with the exercise of computing the
fundamental group of the exterior of a link given by a knot $K$
and its $(p,q)$-cable, and then quotienting by the relations that
$\lambda(K)$ and $\mu(K)$ are trivial.

In principle, we could use the existence of the fibration of
$S^{3} \setminus \nu (K \cup K_{p,q})$, but this requires an
explicit knowledge of the monodromy of the fibration. Another
approach is much more viable. First, it is useful to keep in mind
the construction of $S^{3} \setminus \nu (K \cup K_{p,q})$
obtained by splicing $S^{3} \setminus \hnu K$ and $S^{3} \setminus
\nu (H \cup C_{p,q})$. Now we can proceed (mimicking one of the
standard computations of the fundamental group of a torus knot
exterior) as follows. Presenting the knot and its cable as
companion and satellite, we can write
\begin{equation}
\label{svk} S^{3} \setminus \nu (K \cup K_{p,q}) \;=\; (V_{1}
\setminus \nu K_{p,q}) \:\bigcup_{\Omega}\: (V_{2} \setminus \nu
K_{p,q}).
\end{equation}
The following three paragraphs explain the terms appearing on the right
hand side of (\ref{svk}). 

First, $V_{1}$ is a $3$-manifold with torus boundary that is given
by the exterior $S^{3} \setminus \hnu K$ of the knot $K$\/ union a
collar of its boundary, extended in the outward direction.  $K_{p,q}$ lies on the 
boundary of $V_1$.  By carving out
the intersection with a tubular neighborhood of $K_{p,q}$, we
obtain the manifold $(V_{1} \setminus \nu K_{p,q})$ that is a
deformation retract of $V_{1}$.  Recall that we denote by $\mk$ and
$\lk$ the meridian and the longitude of $S^{3} \setminus \hnu K$.
We point out that these should not be confused with the
meridian and the longitude $\mu(K)$ and $\lambda(K)$ of $S^{3}
\setminus \nu K$; in particular, they are not homotopic to  
$\mu(K)$ and $\lambda(K)$ in
$S^{3} \setminus \nu (K \cup K_{p,q})$.

The 3-manifold $V_{2}$ is a solid torus knotted in $S^{3}$ as a
tubular neighborhood of $K$, with a neighborhood of the core ($\nu K$) 
removed.  $K_{p,q}$ lies on the outer boundary of $V_2$.   
By carving out the intersection with a tubular
neighborhood of $K_{p,q}$ we obtain $(V_{2} \setminus \nu
K_{p,q})$ that is, again, a deformation retract of $V_{2}$ and in
particular is homotopy equivalent to a 2-torus. The generators of
$\pi_1(V_{2} \setminus \nu K_{p,q})$ are given
exactly by $\mu(K)$ and $\lambda(K)$. $V_{1}$ and $V_{2}$
intersect in a 2-torus, which contains the curve $K_{p,q}$, and which
has a natural homotopy basis given by $\mk=\mu(K)$ and
$\lk=\lambda(K)$ (the identification being in $\pi_{1}(V_{1} \cap
V_{2})$), so that $K_{p,q}$ is homotopic to $\mk^{q}\lk^{p} =
\mu(K)^{q}\lambda(K)^{p}$.

If we avoided removing the neighborhood of $K_{p,q}$, the
intersection of $V_{1}$ and $V_{2}$ being a 2-torus, Seifert-Van
Kampen theorem would give the obvious identifications
$\mk=\mu(K)$, $\lk=\lambda(K)$ in the resulting 3-manifold (and
decomposition (\ref{svk}) would simply be a redundant presentation
of $S^{3} \setminus \nu K$). Instead, the intersection of $(V_{1}
\setminus \nu K_{p,q})$ and $(V_{2} \setminus \nu K_{p,q})$ is an
annulus $\Omega$, given by the exterior of $K_{p,q}$ in the torus
$V_{1} \cap V_{2}$. Such an annulus deformation-retracts to its core,
which is homotopic in $V_{1} \cap V_{2}$ to a parallel copy of $K_{p,q}$,
whose image in $\pi_{1}(V_{i} \setminus \nu K_{p,q})$
is given, respectively, by $\mk^{q}\lk^{p}$ and by
$\mu(K)^{q}\lambda(K)^{p}$.  Therefore an application of 
Seifert-Van Kampen theorem gives
\begin{equation} \label{befq}
\pi_1 \big(S^3 \setminus \nu (K \cup K_{p,q}) \big) \:\cong\:
\frac{\pi_{1} (V_{1} ) \ast \big(\zz[\mu(K)] \oplus
\zz[\lambda(K)]\big)}{\mk^{q}\lk^{p} = \mu(K)^{q} \lambda(K)^{p}}.
\end{equation}

We can now complete our argument.  When we quotient by the
relations $\lambda(K) = \mu(K) = 1$\/ in (\ref{befq})  we obtain that
\begin{equation*} \pi_1 \big( E(2)_K \setminus
\nu(\Sigma_g \cup T_{p,q}) \big) \:\cong\: \frac{\pi_{1}
\big(V_{1}\big)}{\mk^{q}\lk^{p}=1}
\end{equation*}
and considering that $V_{1}$ deformation-retracts to $S^{3} \setminus
\hnu K$ with $\mk$ and $\lk$ corresponding to the meridian and the
longitude of $K$\/ respectively, we obtain that
\begin{equation*}
 \pi_1 \big( E(2)_K \setminus \nu(\Sigma_g \cup T_{p,q}) \big) \:\cong\:
 \pi_1 (S^{3}_{q/p}(K)) .
\end{equation*}
This finishes the proof of (\ref{pi_1(half)}).

The second step in the proof of (\ref{pi_1(Xi_q complement)})
consists of applying Seifert-Van Kampen theorem to
decomposition (\ref{decomposition:Xi_q complement}) using
(\ref{pi_1(half)}). Note that the fundamental group of $S^{3}
\setminus \hnu K$ is generated by $\pi_{1}(\sig)$ (a free group
$\mathbb{F}_{2g(K)}$ on $2g(K)$ generators $a_{i},b_{i}$,
$i=1,\dots,g(K)$) and the meridian $\mk$, subject to the monodromy
relations. The relation $\mk^{q}\lk^{p} = 1$ adds an extra
relation between the meridian and the generators of
$\pi_{1}(\sig)$, since $\lk = \prod_{i=1}^{g(K)}
[a_{i},b_{i}]$. In sum, we can write (\ref{pi_1(half)}) as
\begin{eqnarray} \label{each}
&& \pi_1 \big( E(2)_K \setminus
\nu(\Sigma_g \cup T_{p,q}) \big)  \;=\:\:
\langle a_{1},\dots,a_{g(K)},b_{1},\dots,b_{g(K)},\mk\mid \nonumber \\
&&\hspace{1cm} f_{\#} a_{i} = \mk^{-1} a_{i} \mk,\: f_{\#} b_{i} =
\mk^{-1} b_{i} \mk , \: \mk^{q}\lk^{p} = 1
\rangle,
\end{eqnarray}
where $f_{\#}$ denotes the monodromy map of the fiber bundle
$S^3\setminus\hnu K\rightarrow \mk$.

Let $Y$\/ denote the intersection of the  two halves in
(\ref{decomposition:Xi_q complement}).  Precisely, $Y$\/ is
diffeomorphic to the cartesian product $\,(\Sigma_g\setminus
\bigsqcup_{j=1}^{q} \mathbb{D}_{x_j})\times S^1$, where
$\bigsqcup_{j=1}^{q} \mathbb{D}_{x_j}=\Sigma_g\cap(\nu T_{p,q})\,$
is the disjoint union of $q$\/ disks centered at the
intersection points $\{x_1,\dots,x_q\}$ of $\Sigma_g$ with
$T_{p,q}$.  Note that $\Sigma_{g} \setminus \bigsqcup_{j=1}^{q}
\mathbb{D}_{x_j} = \Sigma_{K,q} \cup \Sigma '$, where $\Sigma'$ is the
genus $1$ surface with boundary  described in the proof
of Lemma~\ref{lemma:X_1}.  The loops of the
form\/ $\ast \times S^{1}$, which are copies of the meridian of $\Sigma_{g}$
in $D_{K}$, are contractible in each half-component of 
(\ref{decomposition:Xi_q complement}), as in the proof of
Lemma~\ref{lemma:X_1}.

Remember now that the surface $\sig$ is contained in $(\Sigma_g
\setminus \bigsqcup_{j=1}^{q} \mathbb{D}_{x_j})$, as the latter is
built by gluing to $\sig$ an annulus (see
Figure~\ref{fig:doublehopf}) in $S^{3} \setminus \nu (H \cup
C_{p,q})$ with $q$ disks removed, and then capping off with $\Sigma'$.    
When applied to (\ref{decomposition:Xi_q complement}), 
Seifert-Van Kampen theorem gives, as
$\pi_1 (D_{K} \setminus \nu \Xi_{p,q})$, the free
product of the fundamental groups of the two halves, modulo the
identification induced by the image of the generators of
$\pi_{1}(\Sigma_g\setminus \bigsqcup_{j=1}^{q} \mathbb{D}_{x_j})$.

\begin{lemma}\label{lemma:inclusion}
The fundamental group\/ $\pi_1 \big( E(2)_K \setminus \nu(\Sigma_g
\cup T_{p,q}) \big)$ of each half-component of decomposition\/
$(\ref{decomposition:Xi_q complement})$ is generated by the image
of\/ $\pi_{1}(\Sigma_g\setminus
\bigsqcup_{j=1}^{q} \mathbb{D}_{x_j})$ under the homomorphism 
induced by the inclusion map.
\end{lemma}

\begin{proof}
In order to prove this lemma we need to keep track of the
representatives of the generators of the fundamental group of
each half-component. From presentation (\ref{each}) we see that
the generators for both copies of $\pi_1 ( E(2)_K \setminus
\nu(\Sigma_g \cup T_{p,q}) )$ denoted by $a_{i}$ and $b_{i}$
are contained in the image of $\pi_{1}(\Sigma_g\setminus
\bigsqcup_{j=1}^{q} \mathbb{D}_{x_j})$ for obvious reasons, resulting in the
identification of the generators denoted by the same symbols. The lone
remaining generator, for each half-component, is $\mk$.  Now
we can observe that, as an element of $\pi_{1}(S^{3} \setminus (K \cup
K_{p,q}))$, $\mk$ can be written as a word in the generators of this
group, for which we choose the presentation used in the proof of
Lemma \ref{lemma:X_2}.  By construction this word will be the same
for both halves of decomposition
(\ref{decomposition:Xi_q complement}). It follows that $\mk$, as
an element of $\pi_{1}(S^{3} \setminus (K \cup K_{p,q}))$, can be
written as a word in $\mu(K)$, the $a_{i}$'s, the $b_{i}$'s and in
the $c_{j}$'s. (A more careful study of
$\pi_1 (S^{3} \setminus \nu(H \cup C_{p,q}))$ would enable us to identify
this word, but we will not need its explicit form.) When we
quotient the group $\pi_{1}(S^{3} \setminus (K \cup K_{p,q}))$ by
the relations $\mu(K)  = \lambda(K) = 1$ to obtain $\pi_1 (
E(2)_K \setminus \nu(\Sigma_g \cup T_{p,q}) )$, we see that
(the image of) $\mk$ can be written as a  word in the generators
of (the image of) $\pi_{1}(\Sigma_{K,q})$, and {\it a fortiori}\/ of
$\pi_{1}(\Sigma_g\setminus \bigsqcup_{j=1}^{q} \mathbb{D}_{x_j})$.
As the word is the same on each half of the decomposition,
the generators $\mk$ of two halves are identified.
\end{proof}

As a consequence of Lemma~\ref{lemma:inclusion} 
the free product with amalgamation
identifies all the corresponding generators of $\pi_1$ of two halves of 
(\ref{decomposition:Xi_q complement}) and we get the
simple formula
\begin{equation*}
\pi_1(D_K\setminus \nu \Xi_{p,q})\:\cong\:
\pi_1 \big( E(2)_K \setminus \nu(\Sigma_g \cup T_{p,q}) \big)\:\cong\:
\pi_{1}(S^{3}_{q/p}(K)),
\end{equation*}
with the latter group presented (using standard notation)
in (\ref{stand}). This concludes the proof of
Theorem~\ref{thm:fundamental}.

\section*{Acknowledgments}
We thank Steven Boyer and Olivier Collin for their help in
obtaining some results of \cite{ep:fundamental group} that have
been used here. We also thank Tolga Etg\"u and Ronald Fintushel
for stimulating discussions. We are also grateful to the anonymous
referee for her or his useful observations.

\end{document}